\documentclass[oneside,draft,12pt]{amsart}
\usepackage{amssymb, amscd}
\usepackage[all]{xy}

\theoremstyle{plain}
\newtheorem{theorem}{Theorem}[section]
\newtheorem{lemma}[theorem]{Lemma}

\newtheorem*{theorem*}{}
\newtheorem{proposition}[theorem]{Proposition}
\newtheorem{corollary}[theorem]{Corollary}

\theoremstyle{definition}

\newtheorem{definition}[theorem]{Definition}
\newtheorem{question}[theorem]{Question}
\newtheorem{example}[theorem]{Example}
\newtheorem{remark}[theorem]{Remark}

\DeclareMathOperator{\Tor}{Tor}

\DeclareMathOperator{\Res}{res}

\DeclareMathOperator{\dime}{dim}

\begin{document}

\title{On the Cohomology of Central Frattini Extensions}
\vskip .5in
\author{Alejandro Adem and Jonathan Pakianathan}

\address{Mathematics Department\\
         University of Wisconsin\\
         Madison, Wisconsin, 53706} 
\email{adem@math.wisc.edu, pakianat@math.wisc.edu}
\thanks{The first author was partially supported by the NSF and the NSA}
\vskip .5in
\begin{abstract}
In this paper we provide calculations for the mod $p$ cohomology
of certain $p$--groups, using topological methods.
More precisely, we look at $p$-groups $G$ defined as central extensions
$1\to V\to G\to W\to 1$ of elementary abelian groups
such that $G/[G,G]\otimes\mathbb F_p=W$ and
the defining $k$--invariants span the entire image of the
Bockstein. We show that if $p>dim~V-dim~W+1$, then the mod $p$ cohomology
of $G$ can be explicitly computed as an algebra of the form
${\mathcal P}\otimes A$ where ${\mathcal P}$ is a polynomial ring on
$2$-dimensional generators and $A$ is the cohomology of a compact
manifold which in turn can be computed as the homology of a  
Koszul complex. As an application we provide a complete determination
of the mod p cohomology of the
universal central extension $1\to H^2(W, \mathbb F_p)\to U\to W\to 1$
provided $p>\binom{n}{2}+1$, where $n=dim~W$.
\end{abstract}

\maketitle
\vskip .5in
\section{Introduction}
Obtaining a complete computation for the cohomology of a non--abelian
$p$--group can be quite difficult. In fact very few general computations
exist, especially in the case when $p$ is an odd prime.
In this paper we will consider certain central Frattini extensions of the
form
\[
1\to (\mathbb Z/p)^v\to G\to (\mathbb Z/p)^w\to 1
\]
where $p$ is an odd prime. The condition that the extension be Frattini means 
that $G/[G,G]\otimes\mathbb F_p=(\mathbb Z/p)^w$. 

Our main result will be that under appropriate
conditions on the extension class, a {\it complete} calculation of the
mod $p$ cohomology of $G$ can be obtained {\it provided} $p$ is
sufficiently large. Recall that a Frattini 
extension as above is uniquely determined
by a subspace $K\subset H^2((\mathbb Z/p)^w,\mathbb F_p)$, where 
$H^*((\mathbb Z/p)^w,\mathbb F_p)\cong \Lambda (e_1,\dots ,e_w)
\otimes \mathbb F_p [b_1 ,\dots ,b_w]$, with $|e_i|=1$, and
$b_i = \beta e_i$ ($\beta$ is the Bockstein operator) for
$i=1, \dots , w$. We use the notation above to state our main result.

\newpage
\begin{theorem}
Let $G$ denote a finite $p$--group ($p$ an odd prime)
defined as a
central extension
$1\to V\to G\to W\to 1$ of elementary abelian $p$--groups, where
$v=dim~V, w=dim~W=dim~H_1(G,\mathbb Z/p)$. Assume that
the subspace $K\subset H^2(W,\mathbb F_p)$ defining the extension
contains the entire image of the
Bockstein, i.e., $K$
has a basis of the form $\{ b_1 ,\dots ,b_w, q_1 ,\dots , q_{v-w}\}$,
where $q_1, \dots , q_{v-w}\in \Lambda (e_1, \dots , e_w)$.
Then the 
following hold
\begin{itemize}
\item Every element of order $p$ in $G$ is central.

\item $H^*(G,\mathbb F_p) \cong \mathbb F_p[\zeta_1,\dots,\zeta_v] 
\otimes H^*(M,\mathbb F_p)$ 
as $\mathbb F_p$-algebras, where $M$ is the total space 
of a $(v-w)$-torus bundle over a $w$-torus and $\mathbb F_p[\zeta_1,\dots,
\zeta_v]$ is a polynomial algebra on generators $\zeta_i$ of degree two. 

\item If $p>v-w +1$, then the Lyndon-Hochschild-Serre
spectral sequence for the extension above collapses at $E_3$ and 
\[
H^*(G,\mathbb F_p)/(\zeta_1,\dots ,\zeta_v)\cong 
\Tor_{\mathbb F_p [c_1,\dots ,c_{v-w}]}(\Lambda (e_1,\dots ,e_w),\mathbb F_p)
\]
where $|c_i|=2$, $i=1,\dots ,v-w$
and the $\Tor$ term is determined
by the $k$--invariants, 
$c_i\mapsto q_i\in\Lambda (e_1,\dots ,e_w).$
\end{itemize}
\end{theorem}

The $\Tor$ term can be explicitly
understood as the homology of a Koszul
complex, described as follows. Let $x_1, \dots ,x_{v-w}$ denote
one--dimensional exterior classes, and denote by $K$ the complex
$\Lambda (e_1,\dots ,e_w)\otimes \Lambda (x_1,\dots ,x_{v-w})$,
with differential determined by $\delta (x_i)=q_i$ for
$i=1,\dots ,v-w$. Then the $\Tor$ term above can simply be identified
with the homology of $K$. 

It would seem unlikely that a purely algebraic approach to group
cohomology could have easily led to this result. 
Indeed one of the main drawbacks of general algebraic methods is the lack of
geometric input required to 
determine ambiguities arising from specific differentials in 
the available spectral sequences (a basic computational device).
In contrast, this type of
problem has been considered by topologists in other situations
pertaining to cohomology computations for bundles.
In particular our 
approach provides a geometric alternative to the hypercohomology
spectral sequence developed by Benson and Carlson (see \cite{BC}) in
a very special situation. Our basic contribution is the observation
that for certain group extensions $G$ (such as those above) we can
obtain a geometric model for the basis of $H^*(G,\mathbb F_p)$ as
a free module over a polynomial subring of maximal rank. 

This in turn allows us to introduce effective
techniques from rational homotopy theory and from there outline
conditions which imply the viability of
a complete calculation (via spectral
sequences).
It would of course be interesting to recover the results here using
purely algebraic methods, expressible perhaps in terms of minimal
resolutions. 

A motivation for this paper was the computation of the mod $p$
cohomology of the {\it universal} central extension $U(n,p)$,

\[
1\to H^2((\mathbb Z/p)^n,\mathbb Z/p)\to U(n,p)\to (\mathbb Z/p)^n
\to 1.
\]
Applying our methods 
we obtain

\begin{theorem}
Let $U(n,p)$ denote the universal central extension

\[
1\to (\mathbb Z/p)^{\binom{n+1}{2}}\to U(n,p)\to (\mathbb Z/p)^n\to 1.
\]
If $p>\binom{n}{2}+1$, there is an exact sequence 

\[
0\to (\zeta_1, \dots ,\zeta_{\binom{n+1}{2}})\to H^*(U(n,p),\mathbb F_p)
\to 
\Tor_{\mathbb F_p [c_{ij}]}
(\Lambda (e_1,\dots ,e_n),\mathbb F_p)\to 0
\]
where the $\Tor$ term is determined by $c_{ij}\mapsto e_ie_j$
for $i<j$, $i,j =1,\dots , n$.
\end{theorem}

Using an existing combinatorial computation, we can in fact
make this explicit for most primes (see corollary~\ref{c:u(n,p)}).
The groups $U(n,p)$ seem like basic objects in finite group cohomology.
Further
information on these examples would be very interesting, in particular
it seems plausible to expect that the spectral sequence associated
to the defining extension should collapse at $E_3$ for {\it all primes}.

The computations described here can be expressed as the collapse
at $E_2$ of the Eilenberg-Moore
spectral sequence associated to certain central extensions.
However, in our situation this collapse is equivalent to a collapse
at $E_3$ of the more familiar Lyndon--Hochschild--Serre spectral
sequence (we will often abbreviate this as LHS); 
to simplify matters we have chosen to use this description.

This paper is organized as follows: in \S 2 we provide preliminary
material required in our proofs, in \S 3 we describe 
how the cohomology
calculations can be reduced to computing the cohomology of
a compact manifold, in \S 4 we present the main result using
the complete cohomology information available for large primes. 
In \S 5 we discuss formulas for the Bockstein.

Finally in \S 6 we raise
a general question about the cohomology of central extensions
of elementary abelian $p$--groups.
Throughout this paper $p$ will always denote
an odd prime. 
\vskip .5in

\section{Preliminaries}
In this section we will introduce the basic definitions and background needed
to prove the results in this paper. We refer the reader to \cite{AM} and
\cite{E} for more details.

\begin{definition}
Let 
$1\to\mathbb Z^r\to\Gamma\to\mathbb Z^l\to 1$
be a central extension. We shall say that it is {\it irreducible} if
$H_1(\Gamma,\mathbb Z)=\mathbb Z^l$. Similarly we say that
a central extension $1\to (\mathbb Z/p)^r\to G\to (\mathbb Z/p)^l\to 1$
of elementary abelian $p$--groups
($p$ a prime)
is Frattini if $H_1(G,\mathbb Z/p)=(\mathbb Z/p)^l.$
\end{definition}

Associated to any irreducible central extension as above, we have
a subgroup
\[
K\subset H^2(\Gamma /[\Gamma ,\Gamma ] ,\mathbb Z)
\]
which determines the extension unambiguously up to isomorphism.
If 
$$\kappa_1,\dots ,\kappa_r\in K$$ 
form a $\mathbb Z$--basis, then
we say that they are a complete collection of $k$--invariants
for the extension. A similar convention is made for
central
 Frattini extensions of elementary abelian $p$--groups,
using mod $p$ coefficients.

Let $\pi:\mathbb Z^l\to (\mathbb Z/p^2)^l$ denote the natural
mod $p^2$ quotient map. If
${\mathcal R}\subset H^*((\mathbb Z/p^2)^l,\mathbb Z/p)$
denotes the
subring generated by one--dimensional classes, then $\pi^*$ induces
an isomorphism ${\mathcal R}\to H^*(\mathbb Z^l,\mathbb Z/p)$.

\begin{definition}
Let $K\subset H^2(\Gamma/[\Gamma ,\Gamma ],\mathbb Z)$, $\Gamma$
as above. We define 
$K_p\subset 
H^2(\Gamma /[\Gamma ,\Gamma ]\otimes\mathbb Z/p^2, \mathbb Z/p)$
as $(\pi^*)^{-1}([K])$, where $[K]$ is the mod $p$ reduction of
$K$.
\end{definition}

This defines a central extension
\[1\to (\mathbb Z/p)^r\to G(p)\to (\mathbb Z/p^2)^l\to 1\]
which also fits into an extension of the form
\[
1\to\Gamma (p)\to \Gamma\to G(p)\to 1
\]
derived from the natural projection. 

Central extensions can also be studied geometrically, in particular
the group $\Gamma$ has a classifying space which is the homotopy
fiber of a map
$$\psi: (\mathbb S^1)^l\to (\mathbb CP^{\infty})^r.$$
Let $\rho : (\mathbb S^1)^t\to (\mathbb S^1)^t$ be the $p$--th
power map in each coordinate; it will induce a map
$B\rho :(\mathbb CP^{\infty})^t\to (\mathbb CP^{\infty})^t$.
We can describe the classifying space of $\Gamma (p)$ as the homotopy
fiber of 
$$B\rho\cdot\psi\cdot\rho : (\mathbb S^1)^l\to (\mathbb CP^{\infty})^r.$$
Note that this map is trivial in mod $p$ cohomology, hence the resulting
$k$--invariants are trivial mod $p$.

We examine the $p$--group $G(p)$, but first we recall a definition.

\begin{definition}
A finite group $G$ is said to satisfy the $pC$ condition
if every element
of order $p$ in $G$ is central.
\end{definition}

\begin{proposition}
$G(p)$ is a finite 
$p$--group satisfying the $pC$ condition
which can
be expressed as a Frattini extension
$$1\to (\mathbb Z/p)^{r+l}\to G(p)\to (\mathbb Z/p)^l\to 1.$$
\end{proposition}
\begin{proof}
By construction the $k$--invariants defining $G(p)$ are decomposable
(as sums of products of one--dimensional classes), hence they are elements
in $H^2((\mathbb Z/p^2)^l,\mathbb Z/p)$ which restrict to zero on the
subgroup $(\mathbb Z/p)^l\subset (\mathbb Z/p^2)^l$. Thus the extension
restricted to this subgroup defines an elementary abelian subgroup of
rank equal to $r+l$ which is the kernel of the natural projection
$G(p)\to (\mathbb Z/p)^l$. This of course has maximal rank in $G(p)$
and is central, hence the proposition follows.
\end{proof}

We know by \cite{D} that the groups
$G(p)$ have Cohen--Macaulay mod $p$ cohomology, i.e. $H^*(G(p),\mathbb F_p)$
is free and 
finitely generated over a polynomial subalgebra on $r+l$ generators.
The following lemma (a version of which appears
in \cite{AKM}) explains how to locate a regular sequence; 
we include a proof for the sake of completeness.

\begin{lemma} 
\label{l:varieties}
Let $G$ denote a finite $p$-group
satisfying the pC condition with $E\subseteq G$ the elementary 
abelian subgroup of maximal rank $n$. Let ${\mathcal P}
\subset 
H^*(G,\mathbb F_p)$ be a polynomial subalgebra such that $H^*(E,\mathbb F_p)$
is a finitely generated module over $res^G_E({\mathcal P})$.
Then $H^*(G,\mathbb F_p)$ is a free and finitely generated module over
${\mathcal P}$.
\end{lemma}

\begin{proof} Let ${\mathcal P}= {\mathbb F}_p[\zeta_1,\dots ,\zeta_n]$,
we know that $H^*(G,\mathbb F_p)$ is Cohen-Macaulay
and by a standard result in commutative algebra 
it follows that under that condition the cohomology will be a free module
over any polynomial subring over which it is finitely generated
(see \cite{Serre-loc}). 
Hence 
we only 
need to prove that $H^*(G,\mathbb F_p)$ is a finitely generated $\mathcal
P$-module. To prove this we will use the more geometric language 
of cohomological varieties (see \cite{E} for background).

Let $V_G(\zeta_i)$ denote the homogeneous hypersurface in $V_G$
(the maximal 
ideal spectrum for $H^*(G,\mathbb F_p)$) defined by $\zeta_i$. Then
$H^*(G, \mathbb F_p)$ 
will be finitely generated over $\mathcal P$ if and only if
$V_G(\zeta_1)\cap\dots
\cap V_G(\zeta_n)=\{0\}$. If we represent the class $\zeta_i$ by an
epimorphism $\Omega^{n_i}({\mathbb F}_p)\to {\mathbb F}_p$ 
with kernel $L_{\zeta_i}$, then 
we know that $V_G(\zeta_i)=V_G(L_{\zeta_i})$, the variety associated
to the annihilator 
of $Ext^*_{{\mathbb F}_pG}(L_{\zeta_i},L_{\zeta_i})$. Moreover 
using basic properties of these varieties, we have that
$V_G(\zeta_1)\cap\dots \cap V_G(\zeta_n)=V_G(L_{\zeta_1}\otimes\dots
\otimes L_{\zeta_i})$. Now the cohomological variety of a module will be
$0$ if and only if the module is projective, hence what we need to
prove is that the module
$L_{\zeta_1}\otimes\dots\otimes L_{\zeta_n}$ is projective. However
by Chouinard's Theorem (see \cite{E}) we know that it is enough to 
check this by restricting 
to maximal elementary abelian subgroups; in this case
$E$ is the only such group and projectivity follows from our
hypothesis, as $\Res_E^G({\mathcal P})\subseteq H^*(E,\mathbb F_p)$ is a
polynomial subalgebra  
over which it is finitely generated (note that by Quillen's detection
theorem, the kernel of $\Res^G_E$ is nilpotent, hence $\mathcal P$ embeds in
$H^*(E)$ under this map).
Hence we conclude that $\zeta_1,\dots ,\zeta_n$ form
a homogeneous system of 
parameters and so $H^*(G, \mathbb F_p)$ is free and finitely
generated as a module over $\mathcal P$. 
\end{proof}

\section{Cohomology Calculations}

The goal of this section will be to analyze and in some instances
compute the cohomology ring of the $p$--groups $G(p)$.
We begin by recording the cohomology of $\Gamma (p)$, which
follows directly from the mod $p$ triviality of its defining 
$k$--invariants.

\begin{proposition}
Given $\Gamma (p)$ as above, its mod $p$ cohomology is an
exterior algebra on $r+l$ one--dimensional generators.
\end{proposition}
\begin{proof}
Indeed $\Gamma (p)$ can be expressed as a central extension
$$1\to \mathbb Z^r\to \Gamma (p)\to \mathbb Z^l\to 1$$
which by construction yields a mod $p$ LHS spectral sequence which collapses
at $E_2$.
\end{proof}

\begin{definition}
Given a $pC$ group $G$, let $\Omega_1(G)$ be the maximal 
elementary abelian subgroup of $G$. We will 
let $H^*(\Omega_1(G),\mathbb F_p)_{red}$ be the quotient algebra of 
$H^*(\Omega_1(G),\mathbb F_p)$ modulo the ideal of nilpotent elements. 
This is a graded polynomial algebra on degree 2 generators.
\end{definition}

\begin{definition} 
A $pC$ group $G$ is said to have the $\Omega_1$-extension 
property if the restriction 
$H^2(G;\mathbb F_p) \rightarrow H^2(\Omega_1(G); \mathbb F_p)_{red}$
is onto. 
\end{definition}

The $\Omega_1$-extension property was studied in~\cite{We} 
and \cite{BP}. We will need the following theorem of
T. Weigel which was proven in~\cite{We} using Hopf algebra techniques.
It represents a strengthening of the Cohen-Macaulay property we
previously explained:

\begin{theorem}
\label{thm: omega1}
If $G$ is a $pC$ group with the $\Omega_1$-extension property, 
then we have an isomorphism of $\mathbb F_p$-algebras
$$
H^*(G; \mathbb F_p) \cong 
\mathbb{F}_p[\zeta_1,\dots,\zeta_n] \otimes C^*
$$ 
for some finite dimensional algebra $C^*$. 
Here $n=\dime{\Omega_1(G)}$ and 
the $\zeta_i$ are degree 2 elements which restrict to a generating set
of $H^2(\Omega_1(G), \mathbb F_p)_{red}$.
\end{theorem} 
 
Notice in theorem~\ref{thm: omega1}, that if $\eta_1, \dots, \eta_n$ 
is another regular sequence of degree 2 elements 
which restricts to the generators of 
$H^2(\Omega_1(G), \mathbb F_p)_{red}$, we have that 
$$
(\eta_1,\dots,\eta_n)=(\zeta_1 + c_1, \dots, \zeta_n + c_n)
$$ 
for some $c_i \in C^2$, $1 \leq i \leq n$. 

However by the universal 
property of the tensor product of algebras, we can 
define a graded algebra endomorphism $\Psi$ of $H^*(G, \mathbb{F}_p)$ 
by $\Psi(\zeta_i)=\zeta_i + c_i$ for all $1 \leq i \leq n$ and 
$\Psi|_{C^*}=Id_{C^*}$. Since we can construct an inverse for $\Psi$ 
analogously, $\Psi$ is a graded algebra automorphism of $H^*(G, \mathbb F_p)$
which takes the ideal $(\zeta_1,\dots,\zeta_n)$ to the ideal 
$(\eta_1,\dots,\eta_n)$. Thus we conclude that we have an isomorphism 
of graded algebras
$$
H^*(G,\mathbb F_p)/(\eta_1,\dots,\eta_n) \cong C^*
$$
for any such regular sequence. Furthermore, 
$H^*(G,\mathbb F_p) \cong \mathbb F_p[\eta_1,\dots,\eta_n] \otimes C^*$ 
as algebras.

Our next step is to locate an explicit regular sequence in $G(p)$.

\begin{proposition}
There is a regular sequence of maximal length in $H^*(G(p),\mathbb F_p)$
given by elements $\zeta_1,\dots ,\zeta_{r+l}\in H^2(G(p),\mathbb F_p)$ which 
restrict to generators of $H^2(\Omega_1(G(p)), \mathbb F_p)_{red}$.
Thus $G(p)$ has the $\Omega_1$-extension property.
\end{proposition}
\begin{proof}
We consider the Lyndon--Hochschild--Serre spectral sequence with mod $p$ coefficients
for the group extension 
$$1\to (\mathbb Z/p)^{r+l}\to G(p)\to (\mathbb Z/p)^l\to 1.$$
Recall that $H^*((\mathbb Z/p)^l,\mathbb F_p)\cong \Lambda (e_1,\dots ,e_l)\otimes
\mathbb F_p [b_1,\dots, b_l]$, where, if
$\beta$ denotes the usual Bockstein operator,
$\beta (e_i) = b_i$ for $i=1,\dots ,l$. If
$\gamma :(\mathbb Z/p^2)^l\to (\mathbb Z/p)^l$
is the natural projection, it induces an isomorphism between the subrings
of the mod $p$ cohomologies generated by 1-dimensional 
classes. The $k$--invariants defining the extension $G(p)$
are precisely $(\gamma^*)^{-1}(K_p)\cup\{b_1 ,\dots ,b_l\}$. As a consequence
of this the differential $d_3$ is zero on the two--dimensional polynomial 
generators in the fiber (indeed the ideal generated by the transgressions of
the one--dimensional classes includes the entire ideal generated by the
$b_1,\dots ,b_l$, hence the Bocksteins of these transgressions are in the
ideal already).

Now choose $\zeta_1, \dots,\zeta_{r+l}$ to
be classes in $H^2(G(p),\mathbb F_p)$ restricting to the two--dimensional
permanent cocycle polynomial classes in the edge of the spectral sequence.
These are our desired elements.
\end{proof}

We now assemble the cohomology information above to understand the cohomology
of the finite $p$--groups $G(p)$.

\begin{theorem}
In the mod p LHS spectral sequence associated to 
\[
1\to \Gamma (p)\to \Gamma\to G(p)\to 1
\]
the action of $G(p)$ is homologically trivial
and the one--dimensional generators in the cohomology of $\Gamma (p)$
transgress to a regular sequence $\{ \zeta_i \}_{i=1}^{r+l}$ 
in $H^2(G(p),\mathbb F_p)$,
$E_3=E_{\infty}$, and in particular if $I$ is the ideal 
generated by these transgressions, then  
\[
H^*(G(p),\mathbb F_p)/I
\cong H^*(\Gamma,\mathbb F_p).
\]
Furthermore we have an isomorphism of $\mathbb F_p$-algebras,
$$
H^*(G(p), \mathbb F_p) \cong P^* \otimes H^*(\Gamma, \mathbb F_p).
$$
where $P^*$ is a polynomial algebra on $(r+l)$ degree 2 generators.
\end{theorem}

\begin{proof}
Let $\Gamma_0(p)$ denote the $p$--Frattini subgroup of $\Gamma$, i.e.,
the kernel of the natural projection $\Gamma\to (\mathbb Z/p)^l$.
This group fits into a commutative diagram of
extensions:
\[
\xymatrix{
& & 1 & 1 \\
& & (\mathbb Z/p)^l \ar@{=}[r] \ar[u] &  (\mathbb Z/p)^l \ar[u] \\
1 \ar[r] & \Gamma (p) \ar[r]  & \Gamma \ar[u] \ar[r]	& G(p) \ar[r] \ar[u] &	1 \\
1 \ar[r] & \Gamma (p)\ar[r] \ar@{=}[u] & \Gamma_0(p) \ar[r] \ar[u] & 
(\mathbb Z/p)^{r+l} \ar[r] \ar[u] & 1 \\  
& & 1 \ar[u] & 1 \ar[u] \\
}
\]
Note that $H=(\mathbb Z/p)^{r+l}\subset G(p)$ is $\Omega_1(G(p))$.
As before, the $k$--invariants defining
$\Gamma_0(p)$ are trivial mod $p$, hence 
$H^*(\Gamma_0(p),\mathbb F_p)$ is also an exterior
algebra on $r+l$ one--dimensional generators, and the map
$\phi :\Gamma_0(p)\to H$ induces a surjection
in mod $p$ cohomology. In particular, if we let
\[H^*(H,\mathbb F_p)=\Lambda (e_1,\dots ,e_{r+l})\otimes
\mathbb F_p [\beta e_1,\dots ,\beta e_{r+l}], ~~~~~
H^*(\Gamma_0(p),\mathbb F_p)=\Lambda (y_1,\dots ,y_{r+l})
\]
then we can assume that
$\phi^*(e_i)=y_i$, $\phi^*(\beta e_i)=\beta y_i$.
On the other hand we can choose $q_i\in\Lambda (e_1,\dots ,e_{r+l})$
such that $\phi^*(q_i)=\beta y_i$. Hence we obtain the following
basis for the kernel of $\phi^*$ in dimension equal to two:
$\{b_1-q_1,\dots ,b_{r+l}-q_{r+l}\}$. 

Now consider the
mod $p$ LHS spectral sequence for the bottom row. Since, 
$H^1(H, \mathbb F_p) \cong H^1(\Gamma_0(p),\mathbb F_p)$, the 
five term exact sequence associated to this spectral sequence shows 
that 
$$
H^1(\Gamma(p),\mathbb F_p)^H \overset{d_2}{\cong} 
ker(\phi^*: H^2(H,\mathbb F_p) 
\rightarrow H^2(\Gamma_0(p), \mathbb F_p)).
$$
By comparing dimensions, we see that $H$ must act trivially on 
$H^1(\Gamma(p), \mathbb F_p)$ and hence on $H^*(\Gamma(p), \mathbb F_p)$.
Furthermore, $d_2$ embeds $H^1(\Gamma(p), \mathbb F_p)$ naturally 
as a subspace of $H^2(H, \mathbb F_p)$. 
Since $H$ is central in $G(p)$, 
 we have $G(p)$ acts trivially on 
$H^*(H, \mathbb F_p)$ and we conclude that 
$G(p)$ acts trivially on $H^1(\Gamma(p), \mathbb F_p)$ by the natural 
$d_2$-embedding above, and hence $G(p)$ acts trivially on all of 
$H^*(\Gamma(p), \mathbb F_p)$. 

Let $I^{\prime}$ denote the ideal generated by the $d_2$-transgressions.
Then evidently $I'$ is generated by
the basis of $ker(\phi^*)$ above, which one can easily verify to be a 
regular sequence in the cohomology of $H$. 

Next we consider the spectral sequence for $\Gamma$ given by the middle 
row of our diagram. We have seen that $G(p)$ acts
trivially on $H^*(\Gamma (p), \mathbb F_p)$. Comparing
spectral sequences, we see that $I$ is generated by trangressions
$\zeta_i$, $i=1,\dots ,r+l$ which restrict to a regular sequence
of maximal length in
$H^*(H,\mathbb F_p)$. 

By theorem~\ref{thm: omega1} and the comments immediately following it, 
$H^*(G(p)) \cong P^* \otimes C^*$ 
where $P^*$ is a polynomial algebra on $\zeta_1, \dots, \zeta_{r+l}$.
 From this, it is easy to compute that $E_3=C^*$ and hence $E_3=E_{\infty}$.
Since $E_{\infty}$ is concentrated on a single horizontal row, there 
are no problems in lifting the ring structure to that of 
$H^*(\Gamma, \mathbb F_p)$. 
Hence we conclude that  
$$
H^*(\Gamma, \mathbb F_p) \cong 
H^*(G(p), \mathbb F_p)/(\zeta_1, \dots, \zeta_{r+l}) \cong C^*.
$$
This completes the proof.
\end{proof}

\begin{corollary}
Let $q(t)$ denote the Poincar\'e series for the mod $p$ cohomology
of $\Gamma$ and $p(t)$ the one for the mod $p$ cohomology of $G(p)$.
Then we have
$p(t)=\frac{q(t)}{(1-t^2)^{r+l}}.$
\end{corollary}

In the next section we examine the cohomology of $\Gamma$ using topological
methods.

\section{Cohomology of $\Gamma$ for Large Primes}

As we saw in the previous section, the computation of the 
$\mathbb F_p$-cohomology ring of
$G(p)$ has been reduced to understanding that of $\Gamma$. It turns out
that for sufficiently large primes we have a very precise algebraic
description of the cohomology groups, obtained using bundle theory.
As described in \cite{LP}, an explicit
model for the
classifying space $B\Gamma$ can be easily given; it will be a compact
manifold described as a bundle over
a generalized torus with fibre another generalized torus. This in turn
can be described (up to homotopy) by a map from a generalized torus
to a product of infinite complex projective spaces. The basic calculational
result is summarized in 

\begin{theorem}
Let $\Gamma$ be defined as a central extension of the form
$1\to\mathbb Z^r\to \Gamma\to\mathbb Z^l\to 1$, defined
by a map $\psi :(\mathbb S^1)^l\to (\mathbb CP^{\infty})^r$. Let $p> r+1$ be any 
prime number, then
$$H^*(\Gamma,\mathbb F_p)\cong 
\Tor_{\mathbb F_p [c_1,\dots ,c_r]}(\Lambda (e_1,\dots ,e_l),\mathbb F_p),$$
where the $\Tor$ term is determined
by $\psi^*(c_1),\dots ,\psi^*(c_r)\in\Lambda (e_1,\dots , e_l)$, and 
$\{e_1,\dots ,e_l\}$ are one--dimensional generators for the exterior
algebra $H^*((\mathbb S^1)^l,\mathbb F_p)$
and $\{c_1, \dots ,c_r\}$ are two--dimensional generators
for the polynomial algebra $H^*((\mathbb CP^{\infty})^r,\mathbb F_p)$.
\end{theorem}
\begin{proof}
By a result due to Lambe and Priddy (see \cite{LP}), the LHS spectral
sequence associated to the defining extension above collapses at
$E_3$ if coefficients are taken in $\mathbb Z_{(p)}$ for $p$ sufficiently
large. The improved lower bound $p>r+1$ was obtained in \cite{CP}.
This implies the collapse of the mod $p$ spectral sequence, the statement
is readily derived from the algebraic interpretation of this $E_3$--term.
\end{proof}
\begin{corollary}
If $p$ is sufficiently large, then 
$$H^*(\Gamma ,\mathbb F_p)\cong 
\Tor_{\mathbb Z [c_1,\dots ,c_r]}
(\Lambda_{\mathbb Z}(e_1,\dots ,e_l),\mathbb Z)
\otimes\mathbb F_p.$$
\end{corollary}
\begin{proof}
Indeed, the integral $\Tor$ term can only have $p$--torsion for finitely
many primes $p$, the result follows from the mod $p$ reduction sequence
for $\Tor$.
\end{proof}

\begin{remark} 
In fact \cite{LP} and \cite{CP} prove that the 
$\mathbb Z_{(p)}$--cohomology of $\Gamma$ for $p$ sufficiently large is isomorphic
to the cohomology of the Koszul complex associated to a certain Lie
algebra; this is precisely the $E_3$--term of the LHS spectral sequence.
\end{remark}

We can now state the main result in this paper, which puts together
the different facts we have proved. 

\begin{theorem}
Let $G$ denote a finite $p$--group ($p$ an odd prime)
defined as a
central extension
$1\to V\to G\to W\to 1$ of elementary abelian $p$--groups, where
$v=dim~V, w=dim~W=dim~H_1(G,\mathbb Z/p)$ and $H^*(W,\mathbb F_p)
\cong \Lambda (e_1,\dots ,e_w)\otimes \mathbb F_p [b_1,\dots , b_w]$. 
Assume that
the subspace $K\subset H^2(W,\mathbb F_p)$ defining the extension
contains the entire image of the
Bockstein, i.e., $K$
has a basis of the form $\{ b_1 ,\dots ,b_w, q_1 ,\dots , q_{v-w}\}$,
where $q_1, \dots , q_{v-w}\in \Lambda (e_1, \dots , e_w)$.
Then the 
following hold
\begin{itemize}
\item Every element of order $p$ in $G$ is central.

\item We have an isomorphism of $\mathbb F_p$-algebras: 
$$
H^*(G,\mathbb F_p) \cong \mathbb 
F_p[\zeta_1,\dots,\zeta_v] \otimes H^*(M, \mathbb F_p)
$$
where $M$ is the total space of a $(v-w)$-torus bundle over a $w$-torus
and the degree of the $\zeta_i$ is two.

\item If $p>v-w +1$, then the LHS 
spectral sequence for the extension above collapses at $E_3$ and 
\[
H^*(G,\mathbb F_p)/(\zeta_1,\dots ,\zeta_v)\cong 
\Tor_{\mathbb F_p [c_1,\dots ,c_{v-w}]}(\Lambda (e_1,\dots ,e_w),\mathbb F_p)
\]
where $|c_i|=2$, $i=1,\dots ,v-w$
and the $\Tor$ term is determined
by the $k$--invariants, 
$c_i\mapsto q_i\in\Lambda (e_1,\dots ,e_w).$
\end{itemize}
\end{theorem}

\begin{remark}
Note that it has been shown (see \cite{AK}) that the cohomology of
a $p$--group satisfying the $pC$ condition cannot be detected on proper
subgroups, hence this is an intrinsic calculation.
\end{remark}

\begin{remark}
The $\Tor$ algebra which appears above can be identified
with the mod $p$ cohomology of a compact manifold (as before). Its fundamental
group $\Gamma$ maps onto $G$ inducing a surjection in mod $p$ cohomology.
\end{remark}

We will now exhibit an interesting class of groups to which this theorem
can be applied. In fact the computation of the cohomology of these groups
was a motivation for this work.

\begin{example}
An easy dimension--count shows that 
$$\dime~H^2((\mathbb Z/p)^n,\mathbb Z/p)=\binom{n+1}{2}.$$
Hence we may construct a Frattini 
central extension of the form

$$1\to (\mathbb Z/p)^{\binom{n+1}{2}}\to U(n,p)\to (\mathbb Z/p)^n\to 1$$
which has the following universal property: given any other central 
 Frattini 
extension of the form 
$$1\to V\to U\to (\mathbb Z/p)^n\to 1$$
there exists a central elementary abelian subgroup $E\subset U(n,p)$ such
that $U(n,p)/E\cong U$.
Given this basic property of $U(n,p)$, a natural problem to consider is the
computation of its mod $p$ cohomology. Our results imply
a complete answer for $p$ sufficiently large.
\end{example}

\begin{definition}
Define $H(n)$ to be the
universal central extension 
$$
1\to H^2(\mathbb{Z}^n, \mathbb Z)\to H(n)\to \mathbb{Z}^n\to 1.
$$
$H(n)$ is sometimes called the free two-step nilpotent group on $n$ 
generators.
\end{definition}

\begin{theorem}
\label{thm: Unp}
If $U(n,p)$ denotes the universal central extension

\[
1\to (\mathbb Z/p)^{\binom{n+1}{2}}\to U(n,p)\to (\mathbb Z/p)^n\to 1,
\]
 then there is an isomorphism of graded algebras
$$
H^*(U(n,p), \mathbb F_p) \cong \mathbb F_p[\zeta_1,\dots,
\zeta_{\binom{n+1}{2}}] \otimes H^*(H(n); \mathbb F_p).
$$

\noindent
If $p>\binom{n}{2}+1$, then the cohomology groups above are
determined by the isomorphism 

$$
H^*(H(n), \mathbb{F}_p) \cong \Tor_{\mathbb F_p [c_{ij}]}
(\Lambda (e_1,\dots ,e_n),\mathbb F_p)
$$
where the module structure is specified
by $c_{ij}\mapsto e_ie_j$
for $i<j$, $i,j =1,\dots , n$.

\end{theorem}

It so happens that the corresponding rational $\Tor$ terms have been
computed using representation theory (see \cite{JW}, \cite{Si}). From this
we obtain the following

\begin{corollary}
\label{c:u(n,p)}
Let $U(n,p)$ be the universal central extension. Then, for a fixed
integer $n>0$ and for all but a 
finite number of primes $p$, the Poincar\'e series for 
$H^*(U(n,p),\mathbb F_p)$ is given by $r(t)=\frac{q(t)}{(1-t^2)^m}$
where, if $q(t)=1+a_1t+\dots + a_mt^m$, $m=\binom{n+1}{2}$, then
\[
a_i= \sum_{f+g=i} \sum_{Y_{\lambda}}\prod_{(s,t)\in Y_{\lambda}}
\frac {n+t-s}{h(s,t)}
\]
where $Y_{\lambda}$ ranges over all symmetric, $f+2g$-box,
$f$-hook Young diagrams, and $h(s,t)$ denotes the hooklength
of the box $(s,t)$.
In particular we have
\[
a_1=n, \quad a_2=\frac{n(n+1)(n-1)}{3}, \quad a_3=\frac{n(n^2-1)(3n-4)(n+3)}{60}.
\]
\end{corollary}

Notice how we have obtained a complete calculation for almost every prime.
This type of result would be extremely difficult to observe using traditional
methods in group cohomology, or even computer--assisted calculations. We have
in fact outlined
an {\it effective} method for describing a basis for the cohomology
of an important class of $p$--groups as modules over a polynomial subalgebra
on $2$--dimensional classes.

\section{Bockstein on $H^*(U(n,p), \mathbb F_p)$}
Let $\mathfrak{h}(n)$ be the $\mathbb{F}_p$-Lie algebra with basis elements 
$\{e_k, e_{i,j}\}$ where the indices $i,j,k$ range over the set 
$\{1,\dots,n\}$ and we insist that $i<j$. 
 
The Lie bracket is given by 
$$
[e_i,e_j]=e_{i,j} \text{ and } e_{i,j} \text{ is central}
$$
for $1 \leq i < j \leq n$.
Notice that $\mathfrak{h}(n)$ has dimension $\binom{n+1}{2}$.
Also notice that the $\Tor$ term in Theorem~\ref{thm: Unp} is 
exactly the Lie algebra cohomology of $\mathfrak{h}(n)$. 

As in~\cite{BP} and work of T. Weigel, to every 
$\mathbb F_p$-Lie algebra $\mathfrak{L}$ there 
corresponds a $p$-group of exponent $p^2$ satisfying the $pC$ condition. 

To construct this group, one 
takes a free $\mathbb Z/p^2$-module $K$ of rank equal to the dimension 
of $\mathfrak{L}$. Then one has a canonical surjection 
$\pi: K \to \mathfrak{L}$ and injection $i: \mathfrak{L} \to K$ 
such that $i \circ \pi$ is multiplication by $p$ on $K$.
One defines the group structure on $K$ by
$$
x \cdot y = x + y + i([\pi(x),\pi(y)])
$$
for all $x,y \in K$.

Let $G(\mathfrak{h}(n))$ be the $p$-group associated in this manner 
to $\mathfrak{h}(n)$. Then $G(\mathfrak{h}(n))$ 
fits into a central Frattini extension 
$$
1\to V \to G(\mathfrak{h}(n))\overset{\pi}{\to}  W\to 1
$$
where $V$ and $W$ are elementary abelian $p$-groups of rank $\binom{n+1}{2}$.

Identifying $W$ with $\mathfrak{h}(n)$ 
as in~\cite{BP}, we may look at the subspace 
$S \subset W$ spanned by $\{e_k\}_{1 \leq k \leq n}$. Using this description, 
one can verify easily 
that $\pi^{-1}(S) \cong U(n,p)$. Thus $U(n,p)$ is a subgroup of 
$G(\mathfrak{h}(n))$ which contains $\Omega_1(G(\mathfrak{h}(n)))$.

 From the results in~\cite{BP}, it follows that we 
have an isomorphism of algebras
$$
H^*(G(\mathfrak{h}(n)), \mathbb F_p) \cong \mathbb 
F_p[\zeta_k, \zeta_{i,j}] \otimes 
\wedge^*(x_k, x_{i,j})
$$
where the indices $i,j,k$ range over the set $\{1,\dots,n\}$ and 
we always have $i<j$. 

Furthermore for $p > 3$, the Bockstein is given by the following formulas:
\begin{align*}
\begin{split}
\beta(x_i) = 0 &\text{ and } \beta(x_{i,j}) = -x_ix_j \\
\beta(\zeta_i) = 0 &\text{ and } \beta(\zeta_{i,j}) = \zeta_i x_j - \zeta_j x_i.
\end{split}
\end{align*}
It was also shown that the $\zeta$-elements restrict to a 
regular sequence in the $\mathbb F_p$-cohomology of 
$\Omega_1(G(\mathfrak{h}(n))$.

Now using that $U(n,p)$ is a subgroup of 
$G(\mathfrak{h}(n))$ and that the $\Omega_1$
subgroups of these two groups coincide, by our previous results we have 

$$
H^*(U(n,p), \mathbb F_p) \cong \mathbb F_p[\zeta_k, \zeta_{i,j}] \otimes 
H^*(H(n),\mathbb F_p)
$$
where the polynomial algebra here is the image of the corresponding 
one for $G(\mathfrak{h}(n))$ and we have abused notation and identified the $\zeta$
elements with their restrictions. Notice also that on the 
$H^1(-,\mathbb F_p)$-level, 
the restriction map from $H^*(G(\mathfrak{h}(n)), \mathbb F_p)$ to 
$H^*(U(n,p), \mathbb F_p)$ 
takes $x_{i,j}$ to zero and $\{x_k\}_{k=1}^n$ 
to a basis of $H^1(U(n,p), \mathbb F_p)$. 
Furthermore, the formulas for the Bockstein 
above restrict with the obvious identifications. Thus we obtain:

\begin{proposition}
\label{pro: Bockstein}
For $p>3$, 
$$
H^*(U(n,p),\mathbb F_p) \cong \mathbb F_p[s_k,s_{i,j}] 
\otimes H^*(H(n),\mathbb F_p)
$$
as algebras where as usual the indices $i,j,k$ range over the set 
$\{1,\dots,n\}$ and we always have $i<j$. 
Furthermore there is a basis $\{x_1,\dots,x_n\}$ of 
$H^1(U(n,p),\mathbb F_p)$ such that 
\begin{align*}
\begin{split}
\beta(s_k) = 0 &\text{ and } \beta(s_{i,j}) = s_i x_j - s_j x_i \\
\beta(x_k) = 0 &\text{ and } x_ix_j=0
\end{split}
\end{align*}
for all $1 \leq i < j \leq n$, $1 \leq k \leq n$.

\end{proposition}

\begin{remark}
For large primes $p$, the Bockstein on $H^*(H(n), \mathbb F_p)$ 
vanishes since the integral homology of $H(n)$ 
is finitely generated. However we still 
cannot conclude that proposition~\ref{pro: Bockstein} gives us 
the complete structure of the Bockstein on $H^*(U(n,p), \mathbb F_p)$ since 
the subalgebra of 
$H^*(U(n,p),\mathbb F_p)$ corresponding to 
$H^*(H(n), \mathbb F_p)$ under the isomorphism above is not in general closed
under the higher Bocksteins. 
\end{remark}

\section{Final Remarks and a Problem}
The results in this paper can be thought of as special cases of a more
general collapse theorem for spectral sequences. In fact, given a 
central extension 

$$1\to V\to G\to W\to 1$$
where $V$ and $W$ are both elementary abelian $p$--groups, there is an
Eilenberg--Moore spectral sequence (see \cite{GM})
with $E_2$ term equal to
$\Tor_{H^*(K(V,2),\mathbb F_p)}(H^*(W,\mathbb F_p),\mathbb F_p)$
and converging to $H^*(G,\mathbb F_p)$. Historically these spectral
sequences have been most useful when they collapse at $E_2$. To the
best of our knowledge this occurs for every group extension of this
type whose cohomology has been computed. This leads us to raise the
following somewhat difficult
\vskip .1in
\begin{question} 
Let 
$$1\to V\to G\to W\to 1$$
denote a central extension where both
$V$ and $W$ are elementary abelian $p$--groups. Can the
mod $p$ Eilenberg--Moore spectral sequence fail to collapse at $E_2$?
If so, give reasonable conditions on the $k$--invariants (as above)
which imply
a collapse.
\end{question}

\end{document}